\documentclass[11pt]{article}
\usepackage[cm]{fullpage}
\usepackage{amscd}
\usepackage{amsmath}
\usepackage{amsfonts}
\newcommand{\al}{\alpha}
\newcommand{\be}{\beta}

\newcommand{\ep}{\epsilon}
\newcommand{\mconn}{\overline{\nabla}}
\newcommand{\ddx}[1]{\frac{\partial}{\partial x^{#1}}}

\newcommand{\dd}[1]{\partial_{#1}}

\newcommand{\ddt}[1]{\frac{\partial #1}{\partial t}}

\newcommand{\app}{\approx}
\newcommand{\ov}[1]{\overline{#1}}
\newcommand{\un}[1]{\underline{#1}}
\newcommand{\mg}{\overline{g}}

\begin{document}
\newcounter{remark}
\newcounter{theor}
\setcounter{remark}{0}
\setcounter{theor}{1}
\newtheorem{claim}{Claim}
\newtheorem{theorem}{Theorem}[section]
\newtheorem{proposition}{Proposition}[section]
\newtheorem{lemma}{Lemma}[section]
\newtheorem{defn}{Definition}[theor]
\newtheorem{corollary}{Corollary}[section]
\newenvironment{proof}[1][Proof]{\begin{trivlist}
\item[\hskip \labelsep {\bfseries #1}]}{\end{trivlist}}
\newenvironment{remark}[1][Remark]{\addtocounter{remark}{1} \begin{trivlist}
\item[\hskip
\labelsep {\bfseries #1  \thesection.\theremark}]}{\end{trivlist}}

\centerline{\bf \Large ON DONALDSON'S FLOW OF SURFACES}
\vspace{4pt}
\centerline{\bf \Large  IN A HYPERK\"AHLER FOUR-MANIFOLD}

\bigskip
\begin{center}
\begin{tabular}{ccc}
{\bf Jian Song}  & & {\bf Ben Weinkove\footnotemark} \\
Johns Hopkins University & & Imperial College \\
Department of Mathematics & & Department of Mathematics \\
Baltimore, MD 21218 & & London SW7 2AZ, U.K.
\end{tabular}
\end{center}
\footnotetext{The second author is supported in
part by National Science Foundation grant DMS-05-04285, and is currently on leave from Harvard University, supported by a Royal Society Research Assistantship}
\bigskip

\noindent
{\bf Abstract.} \ We prove some  basic properties of Donaldson's flow of surfaces in a hyperk\"ahler 4-manifold.  When the initial submanifold
is symplectic with respect to one K\"ahler form and Lagrangian with respect
to another,  we show that certain kinds of
 singularities cannot form, and we prove a convergence result under a condition
related to one considered by M.-T. Wang for the mean curvature flow.

\setlength\arraycolsep{2pt}
\addtocounter{section}{1}

\bigskip
\noindent
{\bf 1. Introduction}

\bigskip

In \cite{Do}, Donaldson introduced a number of geometric evolution equations
arising from the framework of moment maps and diffeomorphisms.  If $(S, \rho)$ and
$(M, \omega)$ are two symplectic manifolds then a moment
map can be constructed for
the action of the group of symplectomorphisms of $(S, \rho)$ on the
space of smooth maps $f:S \rightarrow M$.  The gradient flow of the square of the norm of the moment map gives, in many cases, a parabolic flow of maps.
For example, when the two manifolds are
of the same dimension and
K\"ahler this gives rise to a flow of K\"ahler potentials known as the J-flow (see \cite{Ch},
\cite{We}, \cite{SoWe}).

In this paper, we consider the case when $S$ is a Riemann surface with volume form $\rho$ and $(M, \ov{g}, I, J,K)$
 a hyperk\"ahler 4-manifold.    M has three K\"ahler forms $\ov{\omega}_1$, $\ov{\omega}_2$ and $\ov{\omega}_3$, given by
$$\ov{\omega}_1 (X,Y) = \ov{g}(IX, Y), \ \ov{\omega}_2(X,Y) = \ov{g}(JX,Y), \ \ov{\omega}_3 (X,Y) = \ov{g}(KX,Y).$$
Given an immersion $f:S \rightarrow M$, define three functions on $S$ by
\begin{equation} \label{eqnNa}
N_a = \frac{f^*(\ov{\omega}_a)}{\rho}, \quad \textrm{for } a=1, 2, 3.
\end{equation}
Let $\xi_a$, for $a=1,2,3$, be the Hamiltonian vector fields on $S$ associated to $N_a$.  Then Donaldson's flow is  given by
\begin{equation} \label{eqnHflow}
\ddt{f} = I f_* (\xi_1) + J f_* (\xi_2) + K f_* (\xi_3).
\end{equation}
We will call this  the \emph{Hyperk\"ahler mean curvature flow of Donaldson} or, for short, the H-flow.  It is the gradient flow of the functional
$$E(f) = \sum_{a=1}^3 \| N_a \|^2_{L^2(S, \rho)}.$$
Donaldson showed that $f$ is a critical point of $E$ if and only if
its image is a minimal surface and the induced volume form $d\mu$ on
$S$ is a constant multiple of $\rho$.   He also proved that any
compact immersed minimal surface in a hyperk\"ahler 4-manifold with
normal Euler number $e\ge -2$ which is not a complex curve for any
complex structure on $M$ must be not strictly stable with respect to
$E$ (see also the results of \cite{MiWo}).

Donaldson showed that the H-flow is related to the mean curvature flow  of surfaces in a 4-manifold.  Indeed, define a function $\lambda$ on $S$ by
$$\lambda = \frac{d\mu}{\rho}.$$
Then
the H-flow can be written
\begin{equation} \label{eqnHflow1}
\ddt{f} = \lambda \nabla \lambda + \lambda^2 H,
\end{equation}
where $H$ is the mean curvature vector of $S$ in $M$.
This is proved in \cite{Do}. We give another, very explicit, proof
of this (Proposition \ref{prop1}).

The mean curvature flow
of surfaces in a 4-manifold has been studied intensely over
the last several years.  Motivated by a question of Yau on how to deform a symplectic
submanifold into a holomorphic one, Wang \cite{Wa1} showed that if the initial
submanifold is symplectic, it remains so along the
mean curvature flow and that Type I singularities do not form (see also \cite{ChLi1}). Convergence results have been obtained in a number of special cases \cite{Wa1,Wa2},
\cite{TsWa},
\cite{ChLiTi}.  Also, Lagrangian submanifolds remain Lagrangian along the
mean curvature flow \cite{Sm1}.
  This case has been of great interest in light of
the SYZ conjecture \cite{StYaZa} in Mirror Symmetry \cite{ThYa} and there
are a number of promising results  \cite{Sm2, Sm3}, \cite{SmWa},
\cite{Wa3}, \cite{ChLi2}.  However, the Lagrangian mean curvature flow can
exhibit some very bad behavior in general (see e.g. \cite{Ne}).

The problem of finding
necessary and sufficient conditions for the convergence of the mean curvature
flow, even when Lagrangian or symplectic, is still very wide open.
We hope that, as well as being of intrinsic interest, understanding the behavior
of the H-flow may help us to understand  this related problem.

We are mainly interested in the special case of the H-flow when the initial data is Lagrangian with respect to one K\"ahler form and symplectic with respect to another.  The differential form $\theta = \ov{\omega}_2 + i \ov{\omega}_3$ is an $I$-holomorphic (2,0)-form on $M$.  We consider the space $\mathcal{N}$ given by
$$\mathcal{N} = \{ f  : S \rightarrow M \  | \ f \textrm{ is an immersion
and} \ f^*(\theta) = \rho \}.$$
Donaldson \cite{Do} gives a geometric argument to prove the following:

\bigskip
\noindent
{\bf Theorem 1} \  \emph{Let $f_t$ be a solution of the H-flow for $0\le
t \le T$.   If $f_0
\in \mathcal{N}$ then $f_t \in \mathcal{N}$ for $0\le t \le T$.}
\bigskip

In section 3 we present a proof of this using the maximum principle.  If
an H-flow $f_t$ lies in $\mathcal{N}$ we will call it a \emph{Special H-flow}.
This is a flow of Lagrangian submanifolds in $M$.
One might hope that, at least under certain
conditions, the Special H-flow should converge to a parametrised Special Lagrangian submanifold.  These are elements $f$ of $\mathcal{N}$ with $e^{i\alpha}
f^*(\ov{\omega}_1
+ i \ov{\omega}_2) = d\mu$ for some constant angle $\alpha$, or equivalently,
with $\lambda = \textrm{constant}$.  In the case of hyperk\"ahler 4-manifolds, these Special
Lagrangian submanifolds are complex curves with respect to one of
the complex structures.

An interesting property of the Special H-flow is that the induced area form $d\mu$ is bounded from
above and below.  Moreover, we have the following:

\bigskip
\noindent
{\bf Theorem 2} \ \emph{Let $f_t$ be a solution of the Special H-flow on a maximal time interval $[0,T)$, for $0 \le T \le \infty$.  Then
\begin{enumerate}
\item[(i)]  For $t \in [0,T)$, the induced area form $d\mu=d\mu(t)$ on $S$ satisfies
$$\left( \inf_S \lambda|_{t=0}\right) \rho \le d\mu \le \left( \sup_S \lambda|_{t=0}\right)
\rho.$$
\item[(ii)] If there exists a constant $C_0$ such that the norm squared of the second fundamental form $|A|^2$ satisfies
$$\sup_{S} |A|^2 \le C_0, \qquad \textrm{for } t \in [0,T),$$
then there exist constants $C_k$ for $k=1,2, \ldots$ such that
$$\sup_{S} |\nabla^k A|^2 \le C_k, \qquad \textrm{for } t \in [0,T).$$
\item[(iii)] If $T < \infty$ then
$$\lim_{t \rightarrow T} \sup_{S} |A|^2 = \infty.$$
\end{enumerate}}

By analogy with the mean curvature flow \cite{Wa1} we might expect that along
the Special H-flow, Type I singularities do not form.  A Type
I singularity for the H-flow at a time $T$ is a singularity
satisfying the estimate
$$\sup_S |A|^2 \le \frac{C}{T - t}, \quad \textrm{for } t < T,$$
for some constant $C$.  We show that if there is such a singularity then
the function $\lambda$ must exhibit a certain `bad' behavior.  Precisely,
we define a \emph{Type I* singularity} at time $T$  to be a Type I singularity
 with the additional property
that for each $y_0 \in M$ there exist positive constants $\lambda_0$, $\delta$ and $C'$ such that
\begin{equation} \label{eqnholder}
 |\lambda (x) - \lambda_0| \le C'(d(f(x),y_0)^{2\delta} + (T - t)^{\delta})
\end{equation}
for all $(x,t)\in S \times [0,T)$, where $d$ is the distance function
in $(M,g)$.  Then we prove the following result.

\bigskip
\noindent
{\bf Theorem 3} \  \emph{Along the Special H-flow, Type I* singularities do not form.}
\bigskip

We believe that it should be possible to remove the H\"older-type condition
(\ref{eqnholder})
and rule out all Type I singularities.

Finally we prove that in a certain situation, similar to one considered
by Wang \cite{Wa1} for the mean curvature flow, the Special H-flow can indeed converge to a complex curve.

\bigskip
\noindent
{\bf Theorem 4} \ \emph{Suppose that $\ov{\zeta}$ is an anti-self-dual parallel calibrating form on $M$. Then there exists
 $\epsilon>0$ such that if $f_0 \in
\mathcal{N}$ satisfies
$$ \left( \frac{1}{\lambda}\right)|_{t=0}  > 1- \ep \quad \textrm{and} \quad \frac{f_0^* \ov{\zeta}}{d\mu} > 1 -\ep,$$
then the Special H-flow exists for all time and converges in $C^{\infty}$ to a map $f_{\infty} : S \rightarrow M$ such that $f_{\infty}(S)$ is a totally geodesic complex curve in $M$.}

\bigskip

It is not difficult to find simple examples where these initial
conditions are satisfied. Indeed, let $S$ be a torus
$\mathbb{C}/\Lambda$, for some lattice $\Lambda$, and let $\Omega$
be the symplectic form
 on $S$ induced
from the standard one on $\mathbb{C}$.    Let $M = S \times S$ with
the standard metric $\ov{g}$, coordinates $(y^1, \cdots, y^4)$ from
$\mathbb{R}^4$ and the hyperk\"ahler structure defined by
\begin{eqnarray*}
\ov{\omega}_1 & = & dy^1 \wedge dy^2 + dy^3 \wedge dy^4 \\
\ov{\omega}_2 & = & dy^1 \wedge dy^4 + dy^2 \wedge dy^3 \\
\ov{\omega}_3 & = & dy^1 \wedge dy^3 - dy^2 \wedge dy^4.
\end{eqnarray*}
Set $\ov{\zeta} = dy^1 \wedge dy^3 + dy^2 \wedge dy^4$.  Consider a
diffeomorphism $\phi : S \rightarrow S$  preserving $\Omega$ and
define $f_0 : S \rightarrow M$ by $f_0 (x^1, x^2) = (x^1,
\phi(x^1),x^2, \phi(x^2))$. Let $\rho = f_0^* (\ov{\omega}_2)$. Then
$f_0^*(\ov{\omega}_3) =0$ and, if $\phi$ is sufficiently close to
the identity, $f_0$ satisfies the hypotheses of Theorem 4.

The outline of the paper is as follows:  in section 2, some basic properties
of the flow will be presented, including a derivation of (\ref{eqnHflow1});
Theorem 1 is proved in section 3;  a proof of Theorem 2 and the evolution of the second
fundmental form and mean curvature are given in section 4; section 5 contains
a proof of Theorem 3 using a monotonicity formula and a  blow-up argument; finally, Theorem 4 is proved in section 6.

\setcounter{equation}{0}
\setcounter{lemma}{0}
\addtocounter{section}{1}
\bigskip
\bigskip
\pagebreak[3]
\noindent
{\bf 2.  Basic properties of the flow}
\bigskip

Before we derive formulas for this flow, we will describe some notation.  Fix a point $p \in S$.  We will write $g$ for the induced metric $f^* \ov{g}$ on $S$.
 Let $\{ x^1, x^2 \}$ be a normal coordinate system for $S$ at $p$ and let $\{ y^1, \cdots, y^4 \}$ be a normal coordinate system for $M$ at $f(p)$.  We assume that the coordinate system $\{ \dd{1}, \dd{2} \}$ at $p$ is oriented so that
$$\rho (\dd{1}, \dd{2})>0.$$
We will use lower case letters $i, j, \ldots$ to denote the indices 1 and 2, Greek letters $\al, \be, \ldots$ for the indices 3 and 4, and upper case letters $A, B, \ldots$ for the full range of indices 1,2,3,4.  We will denote $\ddx{i}$ by $\dd{i}$ and frequently identify it with $f_* \dd{i}$.  We assume that the $\{ y^A \}$ are chosen so that  $\frac{\partial}{\partial y^i} = f_* \dd{i}$ at $p$, and hence $\frac{\partial f^A}{\partial x^i} = \delta^A_i$ at $p$.
 We will also pick a local orthonormal frame  $\{ e_3, e_4 \}$ for the normal bundle.

Let $\mconn$ and $\nabla$ denote the Levi-Civita connections for $\ov{g}$ on $M$ and $g$ on $S$ respectively.  We use the superscripts $T$ and $N$
to denote the orthogonal projection of a vector in $TM$ to the tangent and
normal bundles of $S$ respectively.
  At $p$, $(\mconn_{\dd{i}} \dd{j} )^T= \nabla_{\dd{i}} \dd{j} = 0$.
The second fundamental form $A : \, TS \times TS \rightarrow NS$ is defined by
$A(X,Y) = (\mconn_X Y)^N.$
Write
$$A(\dd{i}, \dd{j}) = h_{\al ij} e_{\al}, \quad \textrm{for} \quad h_{\al
ij} = \ov{g}( e_{\al}, \mconn_{\dd{i}} \dd{j} ).$$
The mean curvature vector is given at $p$ by $H=H_{\al} e_{\al}$, where
$H_{\al}= g^{ij} h_{\al ij}.$

Recall now that the Hamiltonian vector fields
 $\xi_a$, for $a=1,2,3$, are defined by
\begin{equation} \label{eqnxi}
\rho ( \xi_a, \cdot) = dN_a,
\end{equation}
with the $N_a$ defined by (\ref{eqnNa}).
Note that, calculating at $p$,
\begin{eqnarray*}
N_1^2 + N_2^2 +N_3^2 & = & \lambda^2 \left( (I_1^{\ 2})^2 + (J_1^{\ 2})^2
+ (K_1^{\ 2})^2 \right)
=  \lambda^2.
\end{eqnarray*}
Writing $\xi_a = \xi_a^i \dd{i}$ we have from (\ref{eqnxi}),
$$\rho_{12} \xi_a^1 = \dd{2} N_a, \quad \rho_{12} \xi_a^2 = - \dd{1} N_a,$$
and hence, at $p$,
$$\xi_a = \lambda  (\dd{2} N_a) \dd{1} - \lambda (\dd{1} N_a) \dd{2}, \quad \textrm{for } a=1,2,3.$$

\begin{proposition} \label{prop1}
The H-flow (\ref{eqnHflow}) can be written
\begin{equation} \label{eqnHflow2}
\ddt{f} = \lambda \nabla \lambda + \lambda^2 H.
\end{equation}
\end{proposition}
\begin{proof}
Calculate at a point $p$,
\begin{eqnarray} \nonumber
I f_* \xi_1 & = & \lambda (\dd{2} N_1) I(\dd{1}) - \lambda (\dd{1} N_1)I(\dd{2}) \\ \nonumber
& = & \lambda \dd{2} (\lambda \ov{g}(I (\dd{1}), \dd{2})) I( \dd{1}) - \lambda \dd{1} (\lambda \ov{g}(I (\dd{1}), \dd{2})) I(\dd{2}) \\ \nonumber
& = & \lambda (\dd{2} \lambda) \ov{g}(\dd{2}, I (\dd{1})) I(\dd{1}) + \lambda (\dd{1} \lambda) \ov{g}(\dd{1}, I(\dd{2})) I(\dd{2}) \\ \nonumber
&& \mbox{} - \lambda^2 \ov{g}(\mconn_{\dd{2}} \dd{1}, I(\dd{2})) I(\dd{1}) - \lambda^2 \ov{g}( \mconn_{\dd{1}} \dd{2}, I(\dd{1})) I(\dd{2})\\ \nonumber
&& \mbox{} + \lambda^2 \ov{g}(\mconn_{\dd{1}} \dd{1}, I(\dd{2})) I(\dd{2}) + \lambda^2 \ov{g}( \mconn_{\dd{2}} \dd{2}, I(\dd{1}) ) I(\dd{1})  .
\end{eqnarray}
The same formula holds for $Jf_* \xi_2$ and for $Kf_* \xi_3$, replacing $I$ wherever it occurs with $J$ and $K$ respectively.
Now observe that for $i=1$ or $2$, $\{ \partial_i, I(\partial_i), J(\partial_i), K(\partial_i) \}$ is an orthonormal basis for $T_pM$.  Let $X$ be any  vector in $N_pS$.  A short calculation shows that if we set
\begin{eqnarray*}
Y & = & \ov{g} (X, I(\partial_2)) I (\partial_1) + \ov{g} ( X, I(\partial_1)) I(\partial_2) \\
&& \mbox{} + \ov{g} (X, J(\partial_2)) J (\partial_1) + \ov{g} ( X, J(\partial_1)) J(\partial_2) \\
&& \mbox{} + \ov{g} (X, K(\partial_2)) K (\partial_1) + \ov{g} ( X, K(\partial_1)) K(\partial_2),
\end{eqnarray*}
then $Y$ is orthogonal to $\partial_1$, $I(\partial_1)$, $J(\partial_1)$ and $K(\partial_1)$ and hence is zero.  Set $X = \mconn_{\dd{1}} \dd{2} = \mconn_{\dd{2}} \dd{1}$, and we see that the terms in $\ddt{f}$ involving $X$ vanish (notice that in our coordinates $\nabla_{\dd{i}} \dd{j} = 0$ and so $X$ is indeed orthogonal to $T_pS$). The remaining terms give
\begin{eqnarray*}
If_* (\xi_1) + Jf_* (\xi_2) +Kf_* (\xi_3)
& = & \lambda (\dd{2} \lambda) \dd{2} + \lambda (\dd{1} \lambda) \dd{1} + \lambda^2 \mconn_{\dd{1}} \dd{1} + \lambda^2 \mconn_{\dd{2}} \dd{2} \\
& = & \lambda \nabla \lambda + \lambda^2 H,
\end{eqnarray*}
where we are again making use of the fact that $\{ \partial_i, I(\partial_i), J(\partial_i), K(\partial_i) \}$ is an orthornormal basis for $T_pM$. Q.E.D.
\end{proof}

A solution to this parabolic flow always exists for a short time.

\begin{proposition}
Given any smooth initial map $f_0 \in \mathcal{M}$, there exists $T>0$ such that (\ref{eqnHflow2}) admits a unique smooth solution $f_t \in \mathcal{M}$ for $t\in [0,T)$.
\end{proposition}

\begin{proof} It is well known that the mean curvature flow $\ddt{f} = H$ admits a unique smooth solution for a short time.   Indeed, although the mean curvature flow is not strictly parabolic, it becomes so after modifying by a family of diffeomorphisms of $S$ (for a detailed proof, see e.g. \cite{Zh}).  By the same argument, since $\lambda$ is only first order in $f$, one can see that
$$\ddt{f} = \lambda^2 H,$$
admits a unique smooth solution for a short time.  Define another family of diffeomorphisms $\phi_t : S \rightarrow S$ by
$$\ddt{\phi^i_t} = \lambda \nabla^i \lambda,$$
with $\phi_0 = \textrm{id}.$  Then $\tilde{f}(x,t) = f(\phi(x,t),t)$ solves (\ref{eqnHflow2}).  Q.E.D.
\end{proof}

We will now describe the evolution of the metric $g_{ij} = \ov{g}(\dd{i}, \dd{j})$, the volume form $d\mu$ and $\lambda$ along the H-flow.

\begin{proposition} \label{lemmaddtg}
Along the H-flow,
\begin{enumerate}
\item[(i)]
$\displaystyle{\ddt{} g_{kl} = \nabla_k \nabla_l (\lambda^2) - 2 \lambda^2 H_{\alpha} h_{\al kl}}$
\item[(ii)] $\displaystyle{\ddt{} d\mu = \left( \triangle ( \frac{\lambda^2}{2}
) - \lambda^2 |H|^2 \right) d\mu}$
\item[(iii)] $\displaystyle{\ddt{\lambda^2} = \lambda^2 \triangle \lambda^2 - 2(\lambda^2)^2 |H|^2},$
\end{enumerate}
where, in (i), the repeated index $\al$ indicates that we are summing over
an orthonormal frame for the normal bundle.
\end{proposition}

\begin{proof}
We will work in our usual coordinate system at a point $p$.  Then
\begin{eqnarray} \label{eqnddtg}
\ddt{} g(\dd{k}, \dd{l}) & = & \ov{g}( \mconn_{\nabla (\frac{\lambda^2}{2}) + \lambda^2 H} \dd{k}, \dd{l}) + \ov{g}(\dd{k}, \mconn_{\nabla (\frac{\lambda^2}{2}) + \lambda^2 H} \dd{l}).
\end{eqnarray}
Note that by the torsion-free condition for $\mconn$ and the equality
$$[\nabla (\frac{\lambda^2}{2}) + \lambda^2 H, f_* \dd{j}] = f_* [ \ddt{} , \dd{j} ] = 0,$$
we have
\begin{eqnarray} \nonumber
\ov{g}( \mconn_{\nabla (\frac{\lambda^2}{2}) + \lambda^2 H} \dd{k}, \dd{l})  & = & \ov{g}(\mconn_{\dd{k}}(\nabla (\frac{\lambda^2}{2}) + \lambda^2 H), \dd{l}) \\ \nonumber
& = & \ov{g}( \frac{1}{2} \mconn_{\dd{k}} (\dd{i} (\lambda^2) \dd{i} ) + \dd{k} (\lambda^2) H + \lambda^2 \mconn_{\dd{k}} H, \dd{l}) \\ \nonumber
& = & \frac{1}{2} \dd{k} \dd{l} (\lambda^2) - \lambda^2 \ov{g}( H, \mconn_{\dd{k}} \dd{l}) \\ \nonumber
& = &  \frac{1}{2} \dd{k} \dd{l} (\lambda^2) - \lambda^2 H_{\al} h_{\al kl}.
\end{eqnarray}
Substituting into (\ref{eqnddtg}) gives (i).  (ii) and (iii) follow easily
from (i). Q.E.D.
\end{proof}

Notice that by the maximum principle, we immediately have:

\begin{corollary}
Along the H-flow,
$$\lambda(t) \le \sup_S \lambda|_{t=0}.$$
\end{corollary}

\addtocounter{section}{1}
\setcounter{lemma}{0}
\setcounter{proposition}{0}
\setcounter{equation}{0}
\bigskip
\bigskip
\pagebreak[3]
\noindent
{\bf 3. Preservation of the condition $f^*(\theta) = \rho$ along the H-flow}
\bigskip

In this section we will prove Theorem 1.   Suppose that we have a smooth
solution $f_t$ of the H-flow for $t \in [0,T]$ for some $T>0$ and suppose that the initial data is in $\mathcal{N}$.
For $i=1,2,3$, define functions $\eta_i$ on $S$ by
$$\eta_i = \frac{f^*\ov{\omega}_i}{d\mu}.$$
Then the condition $f^*(\theta) = \rho$ is equivalent to the equations
$$\eta_2= 1/\lambda \quad \textrm{and} \quad \eta_3=0.$$
Set
$$Q = \left( \eta_2 - 1/\lambda \right)^2 + \eta_3^2.$$
We will show that there exists a constant $C$ depending on $f_t$ for $t\in
[0,T]$ such that
\begin{equation} \label{eqnevolve1}
\frac{d}{dt} Q \le \lambda^2 \triangle Q + C Q.
\end{equation}
The theorem will then follow, since by the maximum principle, $$Q(t) e^{-Ct} \le Q(0) = 0.$$

We will now prove (\ref{eqnevolve1}).  Let $\{ x^1, x^2 \}$ be a normal coordinate system for $S$ centered at a point $p$ as before and set
$$\nu_1 = \frac{ K \dd{1} - \eta_3 \dd{2}}{ \sqrt{1-\eta_3^2}}, \qquad \nu_2 = \frac{K \dd{2} + \eta_3 \dd{1}}{\sqrt{1-\eta_3^2}}.$$
Then $\{ \dd{1}, \dd{2}, \nu_1, \nu_2 \}$ gives an orthonormal basis for $T_pM$.
Let $\xi_3 = \sqrt{1 - \eta_3^2(p)}$.  Then at the point $p$, in this basis, $K$ is given by
$$K =  \begin{pmatrix} 0 \, & \, \eta_3 \, & \,  \xi_3 \, & \, 0 \\ -\eta_3 & 0 & 0 & \xi_3 \\ -\xi_3 & 0 & 0 & -\eta_3 \\ 0 & -\xi_3 & \eta_3 & 0 \end{pmatrix}.$$
Now set $a = \ov{\omega}_2 (\dd{1}, \nu_1)(p)$ and $\xi_2 = \ov{\omega}_2 (\dd{1}, \nu_2)(p)$ so that
$$\xi_2^2 = 1 - \eta_2^2 - a^2.$$  Then at the point $p$, $J$ is given by
$$J =  \begin{pmatrix} 0 \, & \, \eta_2 \, & \,  a \, & \, \xi_2 \\ -\eta_2 & 0 & -\xi_2 & a \\ -a & \xi_2 & 0 & -\eta_2 \\ -\xi_2 & -a & \eta_2 & 0 \end{pmatrix}.$$
For $i,j=1,2$, set $h_{3ij} = \ov{g}( \nu_1, \mconn_{\dd{i}} \dd{j})$ and $h_{4ij} = \ov{g}( \nu_2, \mconn_{\dd{i}} \dd{j})$.  Let $\ov{R}$ be the
curvature of $\ov{g}$, given by
$$\ov{R}(X,Y)Z = \mconn_X \mconn_Y Z - \mconn_Y \mconn_X Z - \mconn_{[X,Y]}
Z.$$

We have the following lemma.

\begin{lemma} \label{lemmaevolveeta}
Let $\ov{\omega}$ be parallel 2-form on $M$ and write $\eta = \frac{f^*\ov{\omega}}{d\mu}$.  Then $\eta$ evolves by
\begin{eqnarray*}
\ddt{} \eta & = & \lambda^2 \triangle \eta + \lambda^2 \left( \eta |A|^2 - 2 \ov{\omega}(\nu_1, \nu_2) (h_{3k1} h_{4k2} - h_{4k1} h_{3k2}) \right) \\
&& \mbox{} + 2 \lambda (\dd{1} \lambda) \, \ov{\omega} (H, \dd{2}) + 2 \lambda (\dd{2} \lambda) \, \ov{\omega} (\dd{1}, H) \\
&& \mbox{} + \lambda (\dd{k} \lambda) \, \ov{\omega} (\mconn_{\dd{1}} \dd{k}, \dd{2} ) + \lambda (\dd{k} \lambda) \, \ov{\omega} ( \dd{1}, \mconn_{\dd{2}} \dd{k}) \\
&& \mbox{} -  \lambda^2 \left( \ov{\omega} (( \ov{R} (\dd{k}, \dd{1}) \dd{k})^N, \dd{2}) - \ov{\omega} (( \ov{R}(\dd{k}, \dd{2}) \dd{k})^N, \dd{1}) \right),
\end{eqnarray*}
where we are summing $k$ from 1 to 2.
\end{lemma}
\begin{proof}  We will give only a sketch of the proof, since the argument
is similar to that of Proposition 2.1 in \cite{Wa1}.
Calculate
\begin{eqnarray} \nonumber
\ddt{} f^*\ov{\omega} (\dd{1}, \dd{2}) & = & \ov{\omega}( \mconn_{\lambda^2 H + \nabla
(\frac{\lambda^2}{2})} \dd{1}, \dd{2}) + \ov{\omega}(\dd{1}, \mconn_{\lambda^2
H + \nabla (\frac{\lambda^2}{2})} \dd{2}) \\ \nonumber
& = & \lambda^2 \ov{\omega} (\mconn_{\dd{1}} H, \dd{2}) + \lambda^2 \ov{\omega}(\dd{1},
\mconn_{\dd{2}} H) \\ \nonumber && \mbox{} +  2 \lambda (\dd{1} \lambda) \, \ov{\omega} (H, \dd{2})
 + 2 \lambda (\dd{2} \lambda) \, \ov{\omega}(\dd{1}, H) \\ \nonumber
&& \mbox{} + \ov{\omega}(\mconn_{\dd{1}}
(\nabla(\frac{\lambda^2}{2})), \dd{2}) + \ov{\omega}(\dd{1}, \mconn_{\dd{2}}
(\nabla (\frac{\lambda^2}{2}))) \\ \nonumber
& = & \lambda^2 \ov{\omega} (\mconn_{\dd{1}} H, \dd{2}) + \lambda^2 \ov{\omega}(\dd{1},
\mconn_{\dd{2}} H) \\ \nonumber && \mbox{} +  2 \lambda (\dd{1} \lambda) \, \ov{\omega} (H, \dd{2})
 + 2 \lambda (\dd{2} \lambda) \, \ov{\omega}(\dd{1}, H)  + \eta \triangle (\frac{\lambda^2}{2}) \\  \label{eqnevolveomega}
&& \mbox{} + \lambda (\dd{k} \lambda)
 \ov{\omega}(\mconn_{\dd{1}} \dd{k}, \dd{2}) + \lambda (\dd{k} \lambda) \ov{\omega}
 (\dd{1}, \mconn_{\dd{2}} \dd{k}).
\end{eqnarray}
But from \cite{Wa1}, at $p$,
\begin{eqnarray} \nonumber
\triangle \eta & = & \ov{\omega} (\mconn_{\dd{1}} H, \dd{2}) - \ov{\omega}(\mconn_{\dd{2}}
H, \dd{1}) + \eta |H|^2 \\ \nonumber
&& \mbox{} + \ov{\omega} ((\ov{R}(\dd{k}, \dd{1}) \dd{k})^N, \dd{2}) - \ov{\omega}((\ov{R}(\dd{k},
\dd{2})\dd{k})^N, \dd{1}) \\ \label{eqnlaplacianeta}
&& \mbox{} - \eta |A|^2 + 2 \ov{\omega} (\nu_1, \nu_2)(h_{3k1} h_{4k2} -
h_{4k1} h_{3k2}).
\end{eqnarray}
Then the lemma follows from (\ref{eqnevolveomega}), (\ref{eqnlaplacianeta})
and the evolution of the volume form.  Q.E.D.
\end{proof}

We introduce the following notation:  for functions $A$ and $B$, write $A \app B$ if, at the point $p$,
$$|A-B| \le C(Q + | \nabla \eta_3| + | \nabla (\eta_2 - 1/\lambda)|),$$
for some constant $C$ depending on $f_t$.
Then observe that $\eta_3 \app 0 \app a$ and $\xi_2^2 \app 1 - \eta_2^2 \app 1- 1/\lambda^2$.  Also, if $h_{ijk} = \ov{g} (K \dd{i}, \mconn_{\dd{j}} \dd{k})$ then $h_{3jk} \app h_{1jk}$, $h_{4jk} \app h_{2jk}$ and
$$h_{ijk} \app h_{jik} \app h_{jki}.$$

We have the following lemma.

\begin{lemma} \
\begin{enumerate}
\item[(i)] $\displaystyle{\ddt{} \eta_3 \app \lambda^2 \triangle \eta_3}$
\item[(ii)] $\displaystyle{\ddt{} (\eta_2 - 1/\lambda) \app \lambda^2 \triangle (\eta_2 - 1/\lambda)}.$
\end{enumerate}
\end{lemma}
\begin{proof}  First write $H =H_i dx^i$ where $H_i = g^{jk} h_{ijk}$.  Then note that
\begin{equation} \label{eqnnablalambda}
\nabla \lambda \app - \lambda^2 \xi_2 H,
\end{equation}
To see
this, calculate at $p$
\begin{eqnarray*}
\dd{i}\lambda & \app & - \lambda^2 \ov{\omega}_2 (\mconn_{\dd{i}} \dd{1}, \dd{2}) - \lambda^2
\ov{\omega}_2 (\dd{1}, \mconn_{\dd{i}} \dd{2}) \\
& \app & - \lambda^2 \xi_2 h_{3i1} - \lambda^2 \xi_2 h_{4i2} \\
& \app & -\lambda^2 \xi_2 H_i.
\end{eqnarray*}
For (i) we will make use of Lemma \ref{lemmaevolveeta} with $\ov{\omega}= \ov{\omega}_3$.
Observe that
\begin{eqnarray*}
\lefteqn{2 \lambda (\dd{1} \lambda) \ov{\omega}_3 (H, \dd{2}) + 2 \lambda (\dd{2}
\lambda) \ov{\omega}_3 (\dd{1}, H)  } \\
 & \app  & 2  \lambda^3 \xi_2 H_1 \xi_3
H_2 - 2 \lambda^3 \xi_2 H_2 \xi_3 H_1  =  0,
\end{eqnarray*}
and
\begin{eqnarray*}
\lefteqn{\lambda (\dd{k} \lambda) \ov{\omega}_3 (\mconn_{\dd{1}} \dd{k},
\dd{2}) + \lambda (\dd{k} \lambda) \ov{\omega}_3 (\dd{1}, \mconn_{\dd{2}}
\dd{k} ) } \\
& \app & - \lambda (\dd{k} \lambda) h_{41k}\xi_3 + \lambda (\dd{k} \lambda)
h_{32k} \xi_3 \app 0.
\end{eqnarray*}
We need to deal now with the curvature terms. First observe that the Ricci curvature can be approximated by
$$\ov{\textrm{Ric}} (KX, Y) \app \frac{1}{2} \ov{R} (X,Y, e_A, K e_A),$$
where we are summing over the orthonormal basis $e_A$ for $TM$.  Then
\begin{eqnarray*}
\lefteqn{\ov{\omega}_3 ((\ov{R}(\dd{k}, \dd{1}) \dd{k})^N, \dd{2}) - \ov{\omega}_3
((\ov{R} (\dd{k}, \dd{2}) \dd{k})^N, \dd{1})} \\
 & = & - \xi_3 \ov{R} (\nu_2,
\dd{k}, \dd{k}, \dd{1}) + \xi_3 \ov{R} (\nu_1, \dd{k}, \dd{k}, \dd{2})
\\
&\app & - \ov{\textrm{Ric}} (K \dd{1}, \dd{2}) =0,
\end{eqnarray*}
and (i) follows.  For (ii), calculate:
$$|A|^2 + 2 (h_{3k1} h_{4k2} - h_{4k1} h_{3k2}) \app |H|^2,$$
and
$$ 2 \lambda (\dd{1}\lambda) \ov{\omega}_2 (H, \dd{2}) + 2\lambda (\dd{2}
\lambda) \ov{\omega}_2 (\dd{1}, H) \app - 2 \lambda (\lambda^2 -1 ) |H|^2,$$
and
$$ \lambda (\dd{k} \lambda) \ov{\omega}_2 (\mconn_{\dd{1}}\dd{k}, \dd{2})
+ \lambda (\dd{k} \lambda) \ov{\omega}_2 (\dd{1}, \mconn_{\dd{2}} \dd{k})
\app - \lambda (\lambda^2 - 1) |H|^2.$$
Then since the curvature terms can be neglected as in (i), we see that
$$ \ddt{}\eta_2 \app \lambda^2 \triangle \eta_2 + \lambda |H|^2 - 3\lambda (\lambda^2
-1) |H|^2.$$
But a short calculation shows that
$$ \ddt{} \left( \frac{1}{\lambda}\right) = \lambda^2 \triangle \left(\frac{1}{\lambda}
\right) + \lambda|H|^2 - 3 \frac{|\nabla \lambda|^2}{\lambda},$$
and (ii) follows after making use of (\ref{eqnnablalambda}).  Q.E.D.
\end{proof}

The inequality
 (\ref{eqnevolve1}) and hence Theorem 1 now follow easily.  Indeed, writing $\ddt{}\eta_3 = \lambda^2 \triangle \eta_3 + E_1$ and $\ddt{}(\eta_2 - 1/\lambda) = \lambda^2 \triangle (\eta_2 - 1/\lambda) + E_2$ for $E_1 \app 0 \app E_2$ we see that
\begin{eqnarray*}
\ddt{Q} & = & \lambda^2 \triangle Q - 2 \lambda^2 |\nabla (\eta_2 - \frac{1}{\lambda})|^2
- 2 \lambda^2 |\nabla \eta_3|^2 + 2E_1 (\eta_2 -\frac{1}{\lambda}) + 2 E_2
\eta_3 \\
& \le & \lambda^2 \triangle Q + CQ,
\end{eqnarray*}
where we have made use of the inequality $2ab \le (a^2 + b^2)$ for $a, b \in
\mathbb{R}$.  Q.E.D.

\addtocounter{section}{1}
\setcounter{lemma}{0}
\setcounter{proposition}{0}
\setcounter{equation}{0}
\bigskip
\bigskip
\pagebreak[3]
\noindent
{\bf 4. Evolution equations for the Special H-flow}
\bigskip

In this section we will give the proof of Theorem 2.  As before, let $\{ x^1, x^2 \}$ be a normal coordinate system at a point $p$ for $S$ with $\rho( \dd{1}, \dd{2})>0$.  Define normal vectors $\nu_i = K \dd{i}$ for $i=1,2$.  Then $\{ \dd{1}, \dd{2,} , \nu_1, \nu_2 \}$ gives an basis for the tangent bundle of $M$ and is orthonormal at $p$.    Note that at the point $p$, we have $ (\mconn_{\dd{i}} \nu_j)^N=0$.
In this coordinate system, the second fundamental form is given by
$h_{ijk} = \mg ( \nu_i , \mconn_{\dd{j}} \dd{k}),$
and we define the mean curvature 1-form $H_i dx^i$ by
$H_i = g^{jk} h_{ijk}.$
The second fundamental form is fully symmetric:
$$h_{ijk} = h_{jik} = h_{jki}.$$

Recall from section 3 that $\eta_1$ is defined by
$$\eta_1  =  \frac{f^*(\ov{\omega}_1)}{d\mu}
 =  \frac{N_1}{\lambda}$$

We can write down expressions now for the almost complex structures in terms of $\eta_1$.

\begin{lemma} \label{lemmacoord}
At $p$, $I$, $J$ and $K$ are given by
$$I = \begin{pmatrix} 0 \, & \, \eta_1 \, & \, 0 \, & \, \sqrt{1-\eta_1^2} \\ -\eta_1 & 0 &  -\sqrt{1-\eta_1^2} & 0 \\ 0 & \sqrt{1-\eta_1^2}  & 0 & -\eta_1 \\  -\sqrt{1-\eta_1^2} & 0 & \eta_1 & 0 \end{pmatrix},$$
$$J =  \begin{pmatrix} 0 \, & \,  \sqrt{1-\eta_1^2} \, & \, 0 \, & \, -\eta_1  \\  -\sqrt{1-\eta_1^2} & 0 & \eta_1 & 0 \\ 0 & -\eta_1 & 0 &  -\sqrt{1-\eta_1^2} \\ \eta_1 & 0 &  \sqrt{1-\eta_1^2} & 0 \end{pmatrix}$$
and
$$K =  \begin{pmatrix} 0 \, & \, 0 \, & \,  1 \, & \, 0 \\ 0 & 0 & 0 & 1 \\ -1 & 0 & 0 & 0 \\ 0 & -1 & 0 & 0 \end{pmatrix},$$
\end{lemma}

\begin{proof}
It is trivial to see that $K$ takes the above form.  By definition, $I_{1}^{\ 2} = \eta_1$.  Since
$$\rho(\dd{1}, \dd{2}) = (f^* \ov{\omega}_2) (\dd{1}, \dd{2}) = \ov{g}(J \dd{1}, \dd{2}),$$
and $\rho(\dd{1}, \dd{2})>0$, we have that $J_{1}^{\ 2}>0$.  Hence $J_{1}^{\ 2} = \sqrt{1-\eta_1^2}$, and we may write
$$J =  \begin{pmatrix} 0 \, & \, \sqrt{1-\eta_1^2} \,  & \, 0 \, & \, a \eta_1  \\ -\sqrt{1-\eta_1^2} & 0 & b \eta_1 & 0 \\  0 & - b\eta_1 & 0 &  c \sqrt{1-\eta_1^2} \\  - a \eta_1 & 0 & -c \sqrt{1-\eta_1^2} & 0 \end{pmatrix},$$
where $a$, $b$ and $c$ are each either 1 or -1.  Now calculate
\begin{eqnarray*}
I & = & JK \\
& = &  \begin{pmatrix} 0 \, & \, \sqrt{1-\eta_1^2} \,  & \, 0 \, & \, a \eta_1  \\ -\sqrt{1-\eta_1^2} & 0 & b \eta_1 & 0 \\  0 & - b\eta_1 & 0 &  c \sqrt{1-\eta_1^2} \\  - a \eta_1 & 0 & -c \sqrt{1-\eta_1^2} & 0 \end{pmatrix} \begin{pmatrix} 0 \, & \, 0 \, & \,  1 \, & \, 0 \\ 0 & 0 & 0 & 1 \\ -1 & 0 & 0 & 0 \\ 0 & -1 & 0 & 0 \end{pmatrix} \\
& = & \begin{pmatrix} 0 & - a \eta_1 &  0 &  \sqrt{1-\eta_1^2}  \\ - b \eta_1 & 0  & -\sqrt{1-\eta_1^2} & 0 \\ 0 & -c \sqrt{1-\eta_1^2}  & 0 & -b \eta_1 \\  c \sqrt{1-\eta_1^2} & 0 & -a \eta_1 & 0 \end{pmatrix},
\end{eqnarray*}
from which we obtain $a=-1$, $b=1$ and $c=-1$.  This completes the proof.
\end{proof}

A short calculation shows that, at the point $p$,
$$\dd{k} \lambda =  \lambda^2 \eta_1 H_k.$$
Then part (i) of Theorem 2 follows from the evolution of $\lambda$ and the
maximum principle.

We now calculate the evolution of the second fundamental form.  We will use the convention that an underlined index denotes the application of the complex structure $K$, so that, for example:
$$\ov{R}_{lj\un{i}k} = \ov{g} (\ov{R}(\dd{l}, \dd{j})  \dd{k},  K\dd{i}).$$

\begin{lemma} \label{lemmaevolve}
We have the following formulae for the Special H-flow:

\bigskip
\noindent
(i) The second fundamental form evolves by
\begin{eqnarray} \nonumber
\ddt{} h_{ijk} & = & \lambda^2 \left\{  \triangle h_{ijk}
 +  (3\lambda^2-2)\left( H_i H_l h_{jkl} +  H_j H_l h_{kil} + H_k H_l h_{ijl} \right. \right. \\ \nonumber && \mbox{} \left. +  2 H_i H_j H_k \right)
 + \lambda \eta_1 \left(H_{l,i} h_{jkl} +   H_{l,j} h_{kil} + H_{l,k} h_{ijl}  \right. \\ \nonumber
&& \left. \mbox{}+ 2H_i H_{j,k} + 2H_j H_{k,i} + 2  H_k H_{i,j} + H_l h_{ijk,l}  \right) \\ \nonumber
&& \mbox{} - \left( H_m h_{mli} h_{ljk} + H_m h_{mlj} h_{lki} + H_m h_{mlk}
h_{lij} \right) \\ \nonumber
&& \mbox{}  + 2 \ov{R}_{limj} h_{ k l m} + 2 \ov{R}_{ljmk} h_{i lm} + 2 \ov{R}_{l k m i} h_{j l m}  \\
\nonumber
&& \mbox{} - \ov{R}_{l i l m} h_{jkm} -  \ov{R}_{l j lm} h_{kim} -
 \ov{R}_{lklm} h_{ijm}  \\ \nonumber
&& \mbox{} + h_{lmr} (h_{ijl} h_{kmr} + h_{jkl} h_{imr} + h_{kil} h_{jmr})
- 2 h_{ilm} h_{jmr} h_{krl}
 \\ \label{eqnevolveh}
&& \mbox{} \left.  - (\mconn_{\dd{j}} \ov{R})_{k ll \un{i}} - (\mconn_{\dd{l}} \ov{R})_{jlk\un{i}} \right\}.
\end{eqnarray}
\bigskip
\noindent
(ii) The square of the norm of the second fundamental form evolves by
\begin{eqnarray} \nonumber
\ddt{} |A|^2
& = & \lambda^2 \left\{ \triangle |A|^2 - 2 |\nabla A|^2 + 4(3\lambda^2 -
2) H_i H_j H_k h_{ijk} + 12 \lambda \eta_1 H_i H_{j,k} h_{ijk} \right. \\
\nonumber
&& \mbox{} + 2 \lambda \eta_1 H_l h_{ijk,l} h_{ijk} + 12 \ov{R}_{lkmi} h_{jlm}
h_{ijk} - 6 \ov{R}_{lilm} h_{jkm} h_{ijk} \\ \nonumber
&& \mbox{} + 6 h_{lmr} h_{ijl} h_{kmr} h_{ijk} - 4 h_{ilm} h_{jmr} h_{krl}
h_{ijk} -  2 (\mconn_{\dd{j}} \ov{R})_{kll\un{i}} h_{ijk} \\
&& \mbox{} \left. - 2 (\mconn_{\dd{l}}
\ov{R})_{jlk \un{i}} h_{ijk} \right\}.
\label{eqnevolveAsquared}
\end{eqnarray}
\noindent
(iii) The square of the norm of the mean curvature evolves by
\begin{eqnarray} \nonumber
\ddt{} |H|^2
 & = & \lambda^2 \left\{ \triangle |H|^2 - 2 | \nabla H|^2 + 4(3\lambda^2-2)|H|^4
+ 10 \lambda \eta_1 H_i H_j H_{i,j} \right. \\ \nonumber
&& \mbox{} + 4 \lambda \eta_1 |H|^2 H_{i,i} + 2h_{ilm} h_{jlm} H_i H_j -
2\ov{R}_{lklm} H_m H_k \\ \label{eqnevolveHsquared}
&& \mbox{} \left. - 2 H_k (\mconn_{\dd{i}} \ov{R})_{kll\un{i}} - 2 H_k (\mconn_{\dd{l}}
\ov{R})_{ilk\un{i}} \right\}.
\end{eqnarray}
\end{lemma}
\begin{proof} Calculate at $p$,
\begin{eqnarray} \nonumber
\ddt{} h_{i j k} & = & (\nabla (\frac{\lambda^2}{2}) + \lambda^2 H) ( \ov{g}(\mconn_{\dd{j}} \dd{k}, \nu_{i})) \\
\nonumber
& = & \ov{g}( \mconn_{\dd{j}} \mconn_{\nabla (\frac{\lambda^2}{2}) + \lambda^2 H} \dd{k}, \nu_{i}) + \ov{g}(\ov{R}( \nabla (\frac{\lambda^2}{2}) + \lambda^2 H, \dd{j}) \dd{k}, \nu_{i}) \\ \nonumber
&& \mbox{} +  \ov{g}(\mconn_{\dd{j}} \dd{k}, \mconn_{\nabla (\frac{\lambda^2}{2}) + \lambda^2 H} \nu_{i}) \\
\nonumber
 & = & \ov{g}( \mconn_{\dd{j}} \mconn_{\dd{k}} (\nabla (\frac{\lambda^2}{2}) + \lambda^2 H), \nu_i) + \ov{g}( \ov{R}( \nabla (\frac{\lambda^2}{2}), \dd{j}) \dd{k}, \nu_{i}) \\  \label{eqnddth2}
&& \mbox{} + \lambda^2 \ov{g}( \ov{R}(H, \dd{j}) \dd{k}, \nu_{i}) -  \ov{g}(
K\mconn_{\dd{j}} \dd{k}, \mconn_{\dd{i}} ( \nabla (\frac{\lambda^2}{2}) +
\lambda^2 H)   ). \qquad
\end{eqnarray}
Now calculate
\begin{eqnarray} \nonumber
\ov{g}(\mconn_{\dd{j}} \mconn_{\dd{k}} (\nabla (\lambda^2)), \nu_{i})
& = & \ov{g}( \dd{j} \dd{l} (\lambda^2) \mconn_{\dd{k}} \dd{l}, \nu_{i}) + \ov{g}( \dd{k} \dd{l} (\lambda^2) \mconn_{\dd{j}} \dd{l}, \nu_{i}) \\ \nonumber
&& \mbox{} + \ov{g}( \dd{l} (\lambda^2) \mconn_{\dd{j}} \mconn_{\dd{k}} \dd{l}, \nu_{i}) \\ \label{eqn2derivslambda}
& = & \dd{j} \dd{l} (\lambda^2) h_{i k l} + \dd{k} \dd{l} (\lambda^2) h_{i j l} + \dd{l} (\lambda^2) h_{ i k l,j},
\end{eqnarray}
and
\begin{eqnarray} \nonumber
\ov{g}(\mconn_{\dd{j}} \mconn_{\dd{k}} (\lambda^2 H) , \nu_{i})
& = & \dd{j} \dd{k} (\lambda^2) H_{i} + \dd{k} (\lambda^2) H_{i,j} + \dd{j} (\lambda^2) H_{i,k} \\ \label{eqn2derivsH}
&& \mbox{} + \lambda^2 \ov{g}( \mconn_{\dd{j}} \mconn_{\dd{k}} H , \nu_{i}).
\end{eqnarray}
For the last term, we follow the calculation in \cite{Wa1}, p. 332.  First,
\begin{eqnarray} \nonumber
\ov{g}( \mconn_{\dd{j}} \mconn_{\dd{k}} H, \nu_{i}) & = & \ov{g}(\mconn_{\dd{j}} (\mconn_{\dd{k}} H)^T, \nu_{i}) + \ov{g}(\mconn_{\dd{j}} (\mconn_{\dd{k}} H)^N, \nu_{i}).
\end{eqnarray}
We want to calculate the first term on the right hand side of this equation. To begin with, write
$$\ov{g}( (\mconn_{\dd{k}} H)^T, \dd{l}) = \ov{g}( \mconn_{\dd{k}} H, \dd{l})  =  - \ov{g}(H, \mconn_{\dd{k}} \dd{l})  = -H_{m} h_{m k l},$$
and so
$$\ov{g}(\mconn_{\dd{j}} (\mconn_{\dd{k}} H)^T, \nu_{i}) =  - \ov{g}( (\mconn_{\dd{k}} H)^T, \mconn_{\dd{j}} \nu_{i}) =  -H_{m} h_{m k l} h_{i j l}.$$
Hence
\begin{eqnarray} \label{eqnroughlaplacianH}
\ov{g}( \mconn_{\dd{j}} \mconn_{\dd{k}} H, \nu_{i}) & = & H_{i,kj }  -H_{m} h_{m k l} h_{i j l}.
\end{eqnarray}
The last term of (\ref{eqnddth2}) is given by
\begin{eqnarray} \label{eqnb}
-  \ov{g}(K\mconn_{\dd{j}} \dd{k}, \mconn_{\dd{i}} ( \nabla (\frac{\lambda^2}{2}) +\lambda^2 H)) & = & \frac{1}{2} \dd{i} \dd{l} (\lambda^2) h_{jkl} - \lambda^2
H_l h_{jkm} h_{ilm}. \quad
\end{eqnarray}
Inserting   (\ref{eqn2derivslambda}), (\ref{eqn2derivsH}), (\ref{eqnroughlaplacianH}) and (\ref{eqnb}) into (\ref{eqnddth2}), and making use of the Codazzi equation,
we have
\begin{eqnarray} \nonumber
\ddt{} h_{i j k} & = & \lambda^2 H_{i,kj} + \frac{1}{2} ( \dd{j}\dd{l} (\lambda^2)
h_{ikl} + \dd{k} \dd{l} (\lambda^2) h_{ijl} + \dd{i}\dd{l} (\lambda^2) h_{jkl}
) \\ \nonumber
&& \mbox{}  + \dd{j} \dd{k} (\lambda^2)
H_i + \dd{k} (\lambda^2) H_{i,j}  + \dd{j} (\lambda^2) H_{i,k}
 + \frac{1}{2} \dd{l}(\lambda^2) h_{ijk,l} \\
&& \mbox{} - \lambda^2( H_m h_{mlk}
h_{ijl} +  H_m h_{mli} h_{jkl}) + \lambda^2 H_l \ov{R}_{\un{l} j \un{i}
k}. \label{eqnddth}
\end{eqnarray}
Recall that at the point $p$,
\begin{eqnarray} \label{eqndlambda}
\dd{k} (\lambda^2) & = & 2 \lambda^3 \eta_1 H_k.
\end{eqnarray}
A straightforward calculation shows that
\begin{eqnarray}
\dd{j} \dd{l} (\lambda^2) & = & 2 \lambda^2 (3\lambda^2 - 2 ) H_j H_l + 2 \lambda^3 \eta_1 H_{l,j}, \label{eqnddlambda}
\end{eqnarray}
where we are using the fact that $ \eta_1^2 = 1 - 1/\lambda^2$.
In (\ref{eqnddth}), we obtain
\begin{eqnarray} \nonumber
\ddt{} h_{ijk} & = & \lambda^2 \{ H_{i,kj} + (3\lambda^2 - 2) ( H_i H_l h_{jkl}
+ H_j H_l h_{ikl} + H_k H_l h_{ijl} + 2H_i H_j H_k) \\ \nonumber
&& \mbox{} + \lambda \eta_1 (H_{l,i} h_{jkl} + H_{l,j} h_{kil} + H_{l,k}
h_{ijl} + 2 H_i H_{j,k} + 2H_j H_{k,i} + 2 H_k H_{i,j}  \\
&& \mbox{} + H_l h_{ijk,l}) - ( H_m h_{mlk}
h_{ijl} +  H_m h_{mli} h_{jkl}) + H_l \ov{R}_{\un{l} j \un{i} k} \}. \label{eqnddth3}
\end{eqnarray}
We now use the formula for the Laplacian of the second fundamental form:
\begin{eqnarray} \nonumber
\triangle h_{ijk} & = & H_{i, kj} + (\mconn_{\dd{j}} \ov{R})_{k l l \un{i}} + (\mconn_{\dd{l}} \ov{R})_{jlk \un{i}} \\ \nonumber
&& \mbox{} - 2 \ov{R}_{l k m i} h_{j l m} - 2 \ov{R}_{limj} h_{ k l m} - 2 \ov{R}_{ljmk} h_{i lm} + \ov{R}_{\un{l} j \un{i} k} H_{l} \\
\nonumber
&& \mbox{} + \ov{R}_{l i l m} h_{jkm} +  \ov{R}_{l j lm} h_{kim} +
 \ov{R}_{lklm} h_{ijm} - h_{jkl} \ov{R}_{il} \\ \nonumber
&& \mbox{} + h_{ikm} (h_{lmj} H_l - h_{lmr} h_{ljr}) + h_{iml} (h_{rmj} h_{rkl} - h_{rml} h_{rkj}) \\ \nonumber
&& \mbox{} + h_{m kr} ( h_{m l j} h_{ilr} - h_{mlr} h_{ilj}).
\end{eqnarray}
Substituting into (\ref{eqnddth3}) gives (\ref{eqnevolveh}).

We now calculate the evolution of the square of the second fundamental form,
$$|A|^2 = g^{ip} g^{jq} g^{kr} h_{ijk} h_{pqr}.$$
First, at a point $p$ we have, using (\ref{eqnddlambda}),
\begin{eqnarray} \nonumber
\ddt{} g^{ip} & = & - \nabla_i \nabla_p (\lambda^2) + 2 \lambda^2 H_m h_{mip} \\  \label{eqnddtginverse}
& = & - 2 \lambda^2 ( 3\lambda^2 -2 ) H_i H_p - 2 \lambda^3 \eta_1 H_{i,p} +  2 \lambda^2 H_m h_{mip}.
\end{eqnarray}
Then
\begin{eqnarray*}
\ddt{} |A|^2 & = & 2 \left( \ddt{} h_{ijk} \right) h_{ijk} - 6 \lambda^2 \left(  (  3\lambda^2 -2 ) H_i H_p +  \lambda \eta_1 H_{i,p}  -   H_m h_{mip} \right) h_{ijk} h_{pjk}.
\end{eqnarray*}
Making use of (\ref{eqnevolveh}), we obtain (\ref{eqnevolveAsquared}) after
some simplification.

Finally, we calculate the evolution of the square of the mean curvature.   Using (\ref{eqnddtginverse}),
\begin{eqnarray} \nonumber
\ddt{} |H|^2 & = & \ddt{} (g^{kl} g^{ij} h_{ijk} g^{pq} h_{pql}) \\ \nonumber
& = & \left( \ddt{} g^{kl} \right) H_k H_l + 2 \left( \ddt{} g^{ij} \right) h_{ijk} H_k + 2 g^{ij} H_k\ddt{} h_{ijk} \\ \nonumber
& = & - 2\lambda^2 (3 \lambda^2-2) |H|^4  - 2\lambda^3 \eta_1 H_{l,k} H_k H_l + 2 \lambda^2 H_m h_{mkl} H_k H_l \\ \nonumber
&& \mbox{} - 4 \lambda^2 (3 \lambda^2-2) H_i H_j H_k h_{ijk} - 4 \lambda^3 \eta_1 H_{j,i} h_{ijk} H_k \\  \nonumber
&& \mbox{} + 4 \lambda^2 H_m h_{mij} h_{ijk} H_k + 2g^{ij} H_k \ddt{} h_{ijk}.
\end{eqnarray}
Then from (\ref{eqnevolveh}), we obtain (\ref{eqnevolveHsquared}), after some canceling.  This finishes the proof of the lemma.  Q.E.D.

\end{proof}

\noindent
{\bf Proof of Theorem 2} \ From the evolution equations of $|A|^2$, $\lambda$ and $g$, by a similar argument to \cite{Hu1}, \cite{Hu2}, we obtain for constants
$C(k)$,
$$ \ddt{} |\nabla^k A|^2 \le \lambda^2 \triangle | \nabla^k A|^2 - 2 \lambda^2
|\nabla^k A|^2 + C(k) \sum_{a+b+c=k} |\nabla^a A| \, |\nabla^b A| \, |\nabla^c
A| \, | \nabla^k A|,$$
and (ii) follows by an induction argument.  Part (iii) can be proved using
a straightforward modification of the argument of \cite{Hu2}, section 8.
 Q.E.D.

\bigskip
\bigskip
\pagebreak[3]
\addtocounter{section}{1}
\setcounter{lemma}{0}
\setcounter{theorem}{0}
\setcounter{equation}{0}
\noindent
{\bf 5. Singularity formation}

\bigskip

In this section we give the proof of Theorem 3.  Suppose that there
is a Type I* singularity at time $T$.  We will show that this gives
a contradiction using a monotonicity formula. Fix  $y_0 \in M$. Let
$\iota: M \rightarrow \mathbb{R}^N$ be an isometric embedding and
suppose that $y_0$ corresponds to the origin in $\mathbb{R}^N$.
Write $H\in TM/TS$ for the mean curvature of $S$ in $M$ and
$\tilde{H} \in T\mathbb{R}^N/TS$ for the mean curvature of $S$ in
$\mathbb{R}^N$.  Set $E = H - \tilde{H}$.

Let $\Phi: \mathbb{R}^N \times [0,T) \rightarrow \mathbb{R}$ be the
`backwards heat kernel' on $\mathbb{R}^N$ given by
$$\Phi(y,t) = \frac{1}{4 \pi \lambda_0^2 (T - t)} \exp \left(- \frac{|y|^2}{4 \lambda_0^2 (T-t)}\right).$$
Observe that
$$ \ov{\nabla}^{\mathbb{R}^N} \Phi = - \frac{\Phi}{2 \lambda_0^2 (T-t)} y.$$
Restricted to $S$, $\Phi$ evolves by
\begin{eqnarray} \nonumber
\lefteqn{
 \ddt{} \Phi + \lambda_0^2 \, \triangle \Phi } \\ \nonumber
& = & - \Phi \left[ \frac{ |f^{\perp}|^2}{4 \lambda_0^2 (T -t)^2} + \frac{\lambda^2 H \cdot f^{\perp}}{2 \lambda_0^2 (T-t)} +  \frac{\lambda \nabla \lambda \cdot f^T}{2 \lambda_0^2 (T-t)} + \frac{f^{\perp} \cdot  \tilde{H}}{2 (T-t)} \right],
\end{eqnarray}
where $f^T$ and $f^{\perp}$ are the tangential and normal components
to $S$ of the function $f$, considered as a vector in
$\mathbb{R}^N$.  For $s\ge 1$, calculate
\begin{eqnarray} \nonumber
\lefteqn{ \frac{d}{dt} \int_S \lambda^s \Phi d\mu } \\ \nonumber & = & \int_S \bigg\{ s(\lambda_0^2 - \lambda^2) \lambda^{s-1} \nabla \lambda \cdot \nabla \Phi - (s+1)s \lambda^s | \nabla \lambda|^2 \Phi  \\ \nonumber
&& \mbox{} - (s+1) \lambda^{s+2} \Phi |H|^2  - \frac{1}{2}\lambda^s \Phi \left| \frac{f^{\perp}}{2 \lambda_0 (T-t)} +  \lambda_0 (1+\frac{\lambda^2}{\lambda_0^2})H    \right|^2 \\ \label{eqnmono1}
&& \mbox{}  + \frac{\lambda^s \Phi f^{\perp} \cdot E}{2 (T-t)} - \frac{\lambda^s \Phi |f^{\perp}|^2}{8 \lambda_0^2 (T -t)^2} +  \frac{1}{2}  \Phi \lambda^s \lambda_0^2
(1+ \frac{\lambda^2}{\lambda_0^2})^2   |H|^2  \bigg\}d\mu.
\end{eqnarray}
Now pick $s$ large enough so that
\begin{equation} \label{eqnineq1}
(s+1)\lambda^{s+2} > 1+ \frac{1}{2} \lambda^s \lambda_0^2 (1+\frac{\lambda^2}{\lambda_0^2})^2
\end{equation}
To deal with the first term of (\ref{eqnmono1}), we make use of (\ref{eqnholder}) and the definition of $\Phi$ to see that
\begin{eqnarray} \nonumber
\lefteqn{
\int_S s | \lambda_0^2 - \lambda^2| \lambda^{s-1} | \nabla \lambda \cdot \nabla \Phi | d\mu } \\ \nonumber
& \le & \int_S s \lambda^s | \nabla \lambda|^2 \Phi d\mu + C \int_{|f|^2 \le (T-t)^{1-\delta/2}} (\lambda- \lambda_0) \frac{|f^T|^2}{(T-t)^2} \Phi d\mu \\ \nonumber
&& \mbox{} + C \int_{|f|^2 > (T-t)^{1-\delta/2}} \frac{|f^T|^2}{(T-t)^2} \Phi d\mu \\ \label{eqnineq2}
& \le &  \int_S s \lambda^s | \nabla \lambda|^2 \Phi d\mu+  \frac{C}{(T-t)^{1- \delta/4}}.
\end{eqnarray}
Finally, observe that
\begin{equation} \label{eqnineq3}
\frac{\lambda^s \Phi f^{\perp} \cdot E}{2  (T-t)} \le \frac{\lambda^s \Phi |f^{\perp}|^2}{8 \lambda_0^2 (T-t)^2} + C\Phi,
\end{equation}
since $E$ is bounded.   Then combining (\ref{eqnmono1}),(\ref{eqnineq1}), (\ref{eqnineq2}) and (\ref{eqnineq3}) we obtain
\begin{eqnarray*}
\lefteqn{ \frac{d}{dt} \int_S \lambda^s \Phi d\mu } \\
& \le & \frac{C}{(T-t)^{1- \delta/4}} - \int_S \Phi \left\{ |H|^2 +  \frac{1}{2}
\lambda^s \left| \frac{f^{\perp}}{2 \lambda_0 (T-t)} + \lambda_0 (1+\frac{\lambda^2}{\lambda_0^2})
H  \right|^2 \right\} d\mu.
\end{eqnarray*}
It follows that $\lim_{t \rightarrow T} \int_S \lambda^s \Phi d\mu$ exists.
We can now use a blow-up argument.  First, observe that if the flow
has a singularity at time $T$ then there exists a constant $C>0$ such that
\begin{equation} \label{eqnAlowerbound}
\sup_S |A|^2 \ge \frac{1}{C(T-t)},
\end{equation}
 for all $t$ sufficiently close to $T$.  Indeed if we set $w = \sup_S |A|^2$
 then by Lemma \ref{lemmaevolve} we have
 $$\ddt{w} \le C' (w^2 +1),$$
 and hence
 $$\tan^{-1}(w(t')) - \tan^{-1}(w(t)) \le C'(t' -t), \textrm{ for }
 t < t' <T.$$
 Letting $t' \rightarrow T$ we see that
 $$w(t) \ge \tan (\pi/2 - C'(T-t)) \ge \frac{1}{C(T-t)},$$
giving (\ref{eqnAlowerbound}).

Pick a sequence of points and times $(x_k, t_k) \in S \times [0,T)$ with
$t_k \rightarrow T$ and $x_k \rightarrow x \in S$ so that
$$\nu_k^2 = \sup_{x \in S,\ t \in [0,t_k]} |A|(x,t) = |A|^2(x_k,t_k).$$
We then perform a parabolic scaling around the point $y_0 = f(x) \in M \subset
\mathbb{R}^N$ by
dilating space by a factor $\nu^2_k$ and time by a factor $\nu_k$.  The second
fundamental form is bounded and
by the estimates of Theorem 2 (cf. \cite{Hu2}) so are its covariant derivatives.
  The rescaled
flows converge smoothly and have nonzero second
fundamental form by (\ref{eqnAlowerbound}).  But from the
monotonicity formula above  we see (cf. the argument in \cite{Wa1}) that the blow-up limit satisfies $H=0$ and $f^{\perp}=0$ and hence  is a plane.  This is a contradiction.
Q.E.D.

\bigskip
\bigskip
\pagebreak[3]
\addtocounter{section}{1}
\setcounter{lemma}{0}
\setcounter{theorem}{0}
\setcounter{equation}{0}
\noindent
{\bf 6. Convergence of the flow}

\bigskip

Let $\ov{\zeta}$ be a parallel anti-self-dual calibrating form on $M$ so that $\ov{\omega}_2^2$ and $\ov{\zeta}^2$ determine opposite orientations on $M$.  Write $\beta_1 =1/\lambda= f^*(\ov{\omega}_2)/d\mu$ and $\beta_2 = f^*
\ov{\zeta}/d\mu$.
We can pick normal coordinates $\{x^1, x^2 \}$ on $S$ and local sections of the normal bundle $e_3, e_4$ such that, writing $e_i = f_* \dd{i}$ for $i=1,2$, we have at a point $p$:
$$ ( \ov{\omega}_2 (e_A, e_B) ) = \begin{pmatrix} 0 & \beta_1 & \gamma_1 & 0 \\
- \beta_1 & 0 & 0 & - \gamma_1 \\
- \gamma_1 & 0 & 0 & \beta_1 \\
0 & \gamma_1 & -\beta_1 & 0
\end{pmatrix}$$
with  $\beta_1^2 + \gamma_1^2=1$, and
$$
( \ov{\zeta} (e_A, e_B) ) = \begin{pmatrix} 0 & \beta_2 & \gamma_2 & 0 \\
- \beta_2 & 0 & 0 & \gamma_2 \\
- \gamma_2 & 0 & 0 & -\beta_2 \\
0 & -\gamma_2 & \beta_2 & 0
\end{pmatrix},$$
with $\beta_2^2 + \gamma_2^2=1$ (see \cite{Wa1},
Section 3).  By a straightforward calculation, we have the following lemma.

\begin{lemma} \label{lemmaevolvebeta} In the above coordinates, $\beta_1$ and $\beta_2$ evolve by:
\begin{eqnarray} \nonumber
\ddt{}  \beta_1 & = & \lambda^2 \left\{ \triangle \beta_1 + \beta_1( |A|^2  - 2h_{3k1} h_{4k2} + 2h_{41k}h_{3k2}) - \lambda \gamma_1^2 ( 2 (h_{411}+ h_{312}) H_4  \right. \\ \label{eqnevolvenu1}
&&   \mbox{} + 2(h_{421} + h_{322}) H_3 + \sum_{k=1}^2 (h_{41k} + h_{32k})^2  \,   ) \left. \right\},
\end{eqnarray}
and
\begin{eqnarray} \nonumber
\ddt{}  \beta_2 & = & \lambda^2 \left\{ \triangle \beta_2 + \beta_2( |A|^2  + 2h_{3k1} h_{4k2} - 2h_{41k}h_{3k2}) + \lambda \gamma_1 \gamma_2 ( 2 (h_{411}+ h_{312}) H_4 \right. \\ \label{eqnevolvenu2}
&&  \mbox{} - 2(h_{421} + h_{322}) H_3 + \sum_{k=1}^2 h_{41k}^2 - \sum_{k=1}^2 h_{32k}^2 \,   ) \left. \right\},
\end{eqnarray}
where  $h_{\alpha ij} =\ov{g} (e_{\al}, \mconn_{\dd{i}} \dd{j})$.
\end{lemma}

Set $\mu = \beta_1 + \beta_2$.  Initially, $\beta_1>1-\ep$, $\beta_2>1-\ep$
and hence $\gamma_1, \gamma_2$ are of order $\epsilon$.
  It follows from Lemma \ref{lemmaevolvebeta} and a maximum principle argument that if  $\epsilon$ is sufficiently small then the minimum of $\mu$ is nondecreasing along
the flow.
It follows that
\begin{equation} \label{eqnddtmu}
\ddt{} \mu \ge \lambda^2 ( \triangle \mu + c_1(\ep) \mu |A|^2),
\end{equation}
for some constant $c_1(\epsilon)$ which tends to 1 as $\ep$ tends to zero. We also have
$$ |\nabla \mu|^2 \le c_2(\ep) \mu^2 |A|^2,$$
for some constant $c_2(\ep)$ which tends to zero as $\ep$ tends to zero. From Lemma \ref{lemmaevolve},
\begin{equation} \label{eqn2ff}
\ddt{} |A|^2 \le \lambda^2 \triangle |A|^2 + K_1 |A|^4 + K_2,
\end{equation}
for some constants $K_1$ and $K_2$.  Then the argument in \cite{Wa1}, section
7,  implies that $|A|^2$ is bounded along the flow for $\epsilon$ sufficiently
small.  By Theorem 2, the covariant
derivatives of the second fundamental form are bounded and the flow exists for all time.  From (\ref{eqnddtmu}) and (\ref{eqn2ff}) we obtain for some constants $C_1$ and
$C_2$,
$$\int_0^{\infty} \left( \int_S |A|^2 d\mu \right) dt \le C_1 \quad \textrm{and} \quad \frac{d}{dt} \left(\int_S
|A|^2 d\mu \right) \le C_2.$$
Hence
$$\int_{S} |A|^2 d\mu \rightarrow 0.$$
Moreover, the maps $f_t$ converge in $L^1$ and, by an argument similar to that
in \cite{La}, the submanifolds $f(S)$ converge in $C^{\infty}$ as local graphs
to a smooth submanifold with zero second fundamental form.  It follows that
the maps $f_t$ converge smoothly to a map $f_{\infty}$  whose image is a totally geodesic
complex curve in $M$.
 Q.E.D.

\bigskip
\noindent
{\bf Acknowledgements}  The authors are very grateful to:  Mu-Tao
Wang  for his excellent course at Columbia University on minimal submanifolds and for a number of useful discussions; our former advisor D.H. Phong for his continued help and advice; S.-T. Yau for his support and encouragement;
and Simon Donaldson, Mark Haskins, Richard
Thomas, Tom Ilmanen and Bill Minicozzi for some helpful discussions.

\end{document}